\documentclass[a4paper,12pt]{amsart}
\usepackage{amssymb}
\numberwithin{equation}{section}
\newtheorem{theorem}{Theorem}[section]

\newtheorem{lem}[theorem]{Lemma}

\theoremstyle{remark}
\newtheorem{rem}[theorem]{Remark}

\newcommand{\R}{\mathbb{R}}

\author[J.~Benameur]{Jamel Benameur}
\address{Department of Mathematics, College of Science, King Saud University\\
Riyadh 11451, Kingdom of Saudi Arabia}
\email{\sl jbenameur@ksu.edu.sa}
\thanks{This project was supported by King Saud University, Deanship of Scientific research, College of Science, Research Center.}

\title[Long time decay to the Lei-Lin solution of 3D Navier-Stokes equations]
{Long time decay to the Lei-Lin solution of 3D Navier-Stokes equations}
\date{\today}
\begin{document}
\begin{abstract}
In this paper we prove, if $u\in\mathcal C([0,\infty),{\bf {\mathcal X}^{-1}}(\mathbb R^3))$ is global solution of 3D  Navier-Stokes equations, then $\|u(t)\|_{{\bf {\mathcal X}^{-1}}}$ decays to zero as time goes to infinity. Fourier analysis and standard techniques are used.
\end{abstract}

%@@@@@@@@@@@@@@@@@@@@@@@@@@@@@@@@@@@@%@@@@@@@@@@@@@@@@@@@@@@@@@@@@@@@@@@@@%@@@@@@@@@@@@@@


%@@@@@@@@@@@@@@@@@@@@@@@@@@@@@@@@@@@@%@@@@@@@@@@@@@@@@@@@@@@@@@@@@@@@@@@@@%@@@@@@@@@@@@@@
\maketitle
\tableofcontents

%@@@@@@@@@@@@@@@@@@@@@@@@@@@@@@@@@@@@%@@@@@@@@@@@@@@@@@@@@@@@@@@@@@@@@@@@@%@@@@@@@@@@@@@@@

\section{Introduction}
In this paper we deal with the following 3-D incompressible Navier-Stokes equations:
$$ \left\{\begin{array}{l}
  \displaystyle\partial_t
u-  \nu\Delta u+ (u.\nabla) u=-\nabla p,\quad
\mbox{in}\;\; \mathbb{R}^+\times\mathbb R^3\\
\mbox{div}\; u = 0 \quad
\mbox{in}\; \; \mathbb{R}^+\times\mathbb R^3  \\
u _{\mid t=0}=u^0 \quad\mbox{in}\;\;\mathbb R^3 ,
\end{array}\right.
\leqno(NS) $$ where $\nu>0$ is the viscosity of the fluid, and
$u=u(t,x)=(u_1,u_2,u_3)$ and  $p=p(t,x)$ denote respectively the
unknown velocity and the unknown pressure of the fluid at the point
$(t,x)\in\mathbb{R}^+\times\mathbb R^3$. Here, $u^0=(u^0_1,u^0_2,u^0_3)$ is a given
initial velocity. If the condition is fairly regular, one
can express the pressure using the speed. The study of local existence is studied by several researchers, Leray \cite{Le1,Le2}, Kato \cite{KAT}, ect....\\
The global existence of weak solutions goes back to Leray \cite{Le1} and Hopf \cite{HE}. The global well-posedness of strong solutions for small initial data is due to Fujita and Kato \cite{FK} in the critical Sobolev space $\dot H^{1/2}$ also Chemin \cite{Ch1} has prove the case of $\dot H^s$, $s>1/2$, Kato \cite{K1} in the Lebesgue space $L^3$, and Koch and Tataro \cite{KT} in the space ${\bf BMO}^{-1}$ (see also \cite{MC2,CG,Pl1}). It should be noted, in all these works, that the norms in corresponding spaces of the initial data are assumed to be very small, smaller than the viscosity coefficient $\nu$ multiplied by tiny positive constant $c$. For further results and details can consult the book by Cannone \cite{MC1}. In \cite{zhen}, the authors consider a new critical space that is contains in ${\bf BMO}^{-1}$, where they show it is sufficient to assume the norms of initial data are less that the viscosity coefficient $\nu$. Then the space used in \cite{zhen} is the following
$${\bf {\mathcal X}^{-1}}(\mathbb R^3):=\{f\in\mathcal D'(\mathbb R^3),\,\int_{\mathbb R^3}\frac{|\widehat{f}(\xi)|}{|\xi|}d\xi<\infty\},$$
with the norm
$$\|f\|_{\bf {\mathcal X}^{-1}}=\int_{\mathbb R^3}\frac{|\widehat{f}(\xi)|}{|\xi|}d\xi.$$
We will also use the notation, for $i=0,\,1$,
$${\bf {\mathcal X}}^{i}(\mathbb R^3):=\{f\in\mathcal D'(\mathbb R^3),\,\int_{\mathbb R^3}|\xi|^i|\widehat{f}(\xi)|d\xi<\infty\}.$$
For the small initial data, the authors proved the global existence, precisely:
\begin{theorem}\label{zhen}\cite{zhen} Let $u^0\in {\bf {\mathcal X}^{-1}}(\mathbb R^3)$, such that
$\|u^0\|_{\bf {\mathcal X}^{-1}}<\nu$, then there is a unique $u\in \mathcal C(\mathbb R^+,{\bf {\mathcal X}^{-1}})$ such that $\Delta u\in L^1(\mathbb R^+,{\bf {\mathcal X}^{-1}})$. Moreover, for all
$t\geq0$
\begin{equation}\label{je5}\sup_{0\leq t<\infty}\Big(\|u(t)\|_{\bf {\mathcal X}^{-1}}+(\nu-\|u^0\|_{\bf {\mathcal X}^{-1}})\int_0^t\|\nabla u(t)\|_{L^\infty}\Big)\leq  \|u^0\|_{\bf {\mathcal X}^{-1}}.\end{equation}
\end{theorem}
To show this theorem, the authors used a method of regularization of the initial data $u_0^\lambda=\zeta^\lambda*u^0$, in order to use the standard local existence  theory of the Navier-Stokes equations. They obtain uniform estimates in suitable spaces, and pass to the weak limit as $\lambda$ tends towards zero. If we change this method, and by using Fixed Point Theorem on $\mathcal C([0,T],{\bf {\mathcal X}^{-1}}(\mathbb R^3))\cap L^1([0,T],{\bf {\mathcal X}^{1}}(\mathbb R^3))$ and Lemma \ref{scl1}, we can deduce the following: Let $u^0\in {\bf {\mathcal X}^{-1}}(\mathbb R^3)$, such that
$\|u^0\|_{\bf {\mathcal X}^{-1}}<\nu$, then there is a unique $u\in \mathcal C(\mathbb R^+,{\bf {\mathcal X}^{-1}})$ such that $\Delta u\in L^1(\mathbb R^+,{\bf {\mathcal X}^{-1}})$. Moreover, for all
$t\geq0$
\begin{equation}\label{jam1}\|u(t)\|_{\bf {\mathcal X}^{-1}}+(\nu-\|u^0\|_{\bf {\mathcal X}^{-1}})\int_0^t\|u(z)\|_{\bf {\mathcal X}^{1}}dz\leq  \|u^0\|_{\bf {\mathcal X}^{-1}}.\end{equation}
Moreover, in \cite{ZY1} Zhang and Yin prove the local existence for large initial data and blow up criteria if the maximal time is finite, precisely:
\begin{theorem}\label{ZYth} Let $u^0$ be in ${\bf {\mathcal X}^{-1}}(\mathbb R^3)$. There exists time $T$ such that the system $(NS)$ has unique solution $u$ in $L^2([0,T],{\bf {\mathcal X}^{-1}}(\mathbb R^3))$ which also belongs to
 $${\mathcal C}([0,T];{\bf {\mathcal X}^{-1}}(\mathbb R^3))\cap L^1([0,T],{\bf {\mathcal X}^1}(\mathbb R^3))\cap L^\infty([0,T];{\bf {\mathcal X}^{-1}}(\mathbb R^3)).$$
 Let $T^*$ denote the maximal time of existence of such solution. Then:\\
 (i) If $\|u^0\|_{{\bf {\mathcal X}^{-1}}}<\nu$, then
 $$T^*=\infty.$$
 (ii) If $T^*$ is finite, then
 $$\int_0^{T^*}\|u(t)\|_{{\bf {\mathcal X}^0}}^2dt=\infty.$$
\end{theorem}
Our main result is to prove non-blowup at large time and the norm of the global solution in ${\bf {\mathcal X}^{-1}}(\mathbb R^3)$ goes to zero at infinity.
\begin{theorem}\label{jamel1} Let $u\in{\mathcal C}(\mathbb R^+,{\bf {\mathcal X}^{-1}}(\mathbb R^3))$ be a global solution of $(NS)$, then
\begin{equation}\label{jamelF1}\limsup_{t\rightarrow\infty}\|u(t)\|_{\bf {\mathcal X}^{-1}}=0.\end{equation}
\end{theorem}
In the following we give a natural application of Theorem \ref{jamel1}, it is the stability of global solutions of $(NS)$ system.
\begin{theorem}\label{jamel2} Let $u\in{\mathcal C}(\mathbb R^+,{\bf {\mathcal X}^{-1}}(\mathbb R^3))$ be a global solution of $(NS)$, then for all $v^0\in{\bf {\mathcal X}^{-1}}(\mathbb R^3)$ such that
$$\|v^0-u(0)\|_{{\bf {\mathcal X}^{-1}}}<\frac{\nu}{8}e^{-\frac{2}{\nu}\int_0^\infty\|\widehat{u}(s)\|_{L^1}^2ds}.$$
Then, Navier-Stokes system starting by $v^0$ has a global solution. Moreover, if $v$ is the corresponding global solution; then, for all $t\geq0$,
$$\|v(t)-u(t)\|_{{\bf {\mathcal X}^{-1}}}+\frac{\nu}{2}\int_0^t\|v(s)-u(s)\|_{{\bf {\mathcal X}^{1}}}ds\leq \|v^0-u(0)\|_{{\bf {\mathcal X}^{-1}}}e^{\frac{2}{\nu}\int_0^\infty\|\widehat{u}(s)\|_{L^1}^2ds}.$$
\end{theorem}
The remainder of this paper is organized in the following way: In section 2 we give some notations and important preliminaries results. Section 3 is devoted to prove the principle result. In section 4 we prove the stability result for global solutions.

\section{Notations and Preliminaries Results}
\subsection{Notations}
\noindent In this short section we collect some notations and
definitions that will be used later on.\\
$\bullet$ The Fourier transformation is normalized as
$${\mathcal
F}(f)(\xi)=\stackrel{\wedge}{f}(\xi)=\displaystyle\int_{\R^3}\exp(-ix.\xi)f(x)dx,\;\;\;
\xi=(\xi_1,\xi_2,\xi_3)\in\R^3.$$ $\bullet$ The inverse Fourier
formula is
$${\mathcal
F}^{-1}(g)(x)=\displaystyle(2\pi)^{-3}\int_{\R^3}\exp(i\xi.x)f(\xi)d\xi,\;\;\;
x=(x_1,x_2,x_3)\in\R^3.$$ $\bullet$ For $s\in\R$, $H^s(\R^3)$
denotes the usual non homogeneous Sobolev space on $\R^3$ and
$\langle.,.\rangle_{H^s(\R^3)}$ denotes the usual scalar product on
$H^s(\R^3)$.\\
$\bullet$ For $s\in\R$, $\dot H^s(\R^3)$ denotes the usual
homogeneous Sobolev space on $\R^3$ and $\langle.,.\rangle_{\dot H^s(\R^3)}$
denotes the usual scalar product on
$\dot H^s(\R^3)$.\\
$\bullet$ The convolution product of a suitable pair of functions
$f$ and $g$ on $\R^3$ is given by $$(f*g)(x):=\int_{\mathbb R^3}f(y)g(x-y)dy.$$\\
$\bullet$
If $f=(f_1,f_2,f_3)$ and $g=(g_1,g_2,g_3)$ are two vector fields,
we set $$f\otimes g:=(g_1f,g_2f,g_3f),$$ and$${\rm div}(f\otimes
g):=({\rm div}(g_1f),{\rm div}(g_2f),{\rm div}(g_3f)).$$ $\bullet$
For any subset $X$ of a set $E$, the symbol ${\bf 1}_X$ denote the
characteristic function of $X$ defined by $${\bf 1}_X(x)=1\;\;
\mbox{if}\;\; x\in X,\;\;\;{\bf 1}_X(x)=0\; \mbox{elsewhere}.$$
\subsection{Preliminaries Results}
\begin{lem}\label{scl1}
$(i)$ If $f,g\in{\bf {\mathcal X}^0}(\mathbb R^3)$, then $fg\in {\bf {\mathcal X}^0}(\mathbb R^3)$ and
$$\|fg\|_{\bf {\mathcal X}^0}\leq \|f\|_{\bf {\mathcal X}^0}\|g\|_{\bf {\mathcal X}^0}.$$
$(ii)$ If $f\in{\bf {\mathcal X}^{-1}}(\mathbb R^3)\cap {\bf {\mathcal X}^{1}}(\mathbb R^3)$, then $f\in{\bf {\mathcal X}^{-1}}(\mathbb R^3)$ and
$$\|f\|_{\bf {\mathcal X}^0}\leq \|f\|_{\bf {\mathcal X}^{-1}}^{1/2}\|f\|_{\bf {\mathcal X}^1}^{1/2}.$$
\end{lem}
{\bf Proof of lemma \ref{scl1}.} $(i)$ is a given by direct application of Young inequality.\\
To prove $(ii)$, we can write
$$\begin{array}{ccc}
\|f\|_{\bf {\mathcal X}^0}&=&\displaystyle\int_{\xi}|\widehat{f}(\xi)|d\xi\\
&=&\displaystyle\int_{\xi}|\xi|^{1/2}|\widehat{f}(\xi)|^{1/2}\frac{|\widehat{f}(\xi)|^{1/2}}{|\xi|^{1/2}}d\xi.
\end{array}
$$
Cauchy-Schwartz inequality gives the desired result.\\
\begin{lem}\label{scl2} If $s>1/2$, we have $H^s(\mathbb R^3)\hookrightarrow {\bf {\mathcal X}^{-1}}(\mathbb R^3)$ and
$$\|f\|_{\bf {\mathcal X}^{-1}}\leq C_s\|f\|_{L^2}^{1-\frac{1}{2s}}\|f\|_{\dot H^s}^{\frac{1}{2s}}.$$
\end{lem}
{\bf Proof of lemma \ref{scl2}.}\\
For $R>0$, we have
$$\|f\|_{\bf {\mathcal X}^{-1}}\leq \|f{\bf 1}_{\{|D|<R\}}\|_{\bf {\mathcal X}^{-1}}+\|f{\bf 1}_{\{|D|>R\}}\|_{\bf {\mathcal X}^{-1}}.$$
Cauchy-Schwartz inequality gives
$$\begin{array}{ccc}
\|f{\bf 1}_{\{|D|<R\}}\|_{\bf {\mathcal X}^{-1}}&=&\displaystyle\int_{|\xi|<R}\frac{|\widehat{f}(\xi)|}{|\xi|}d\xi\\
&\leq&\displaystyle\Big(\int_{|\xi|<R}\frac{1}{|\xi|^2}d\xi\Big)^{1/2}\|f\|_{L^2}\\
&\leq&\displaystyle\sqrt{4\pi} R^{\frac{1}{2}} \,\|f\|_{\dot H^s},\\
  \end{array}
$$
and
$$\begin{array}{ccc}
\|f{\bf 1}_{\{|D|<R\}}\|_{\bf {\mathcal X}^{-1}}&=&\displaystyle\int_{|\xi|>R}\frac{1}{|\xi|^{s+1}}|\xi|^s|\widehat{f}(\xi)|d\xi\\
&\leq&\displaystyle\Big(\int_{|\xi|>R}\frac{1}{|\xi|^{2s+2}}d\xi\Big)^{1/2}\|f\|_{\dot H^s}\\
&\leq&\displaystyle\sqrt{\frac{4\pi}{2s-1}} R^{\frac{1}{2}-s} \,\|f\|_{L^2},\\
  \end{array}
$$
To conclude, it suffices to take $R=\big(\frac{\|f\|_{\dot H^s}}{\|f\|_{L^2}}\big)^{1/s}$.
\endproof
\begin{rem} In the case $ s =1/2$ there is no comparison between $H^{1/2}(\mathbb R^3)$ and ${\bf {\mathcal X}^{-1}}(\mathbb R^3).$ It suffices to consider the functions $f$ and $g$ defined as follows
$$f=\mathcal F^{-1}\Big(\frac{1}{|\xi|^{3/2}}{\bf 1}_{\{|\xi|<1\}}\Big)\;\;{\rm and }\;\;g=\mathcal F^{-1}\Big(\frac{1}{|\xi|^{7/4}}{\bf 1}_{\{|\xi|>1\}}\Big).$$
Indeed:
$$\|f\|_{\bf {\mathcal X}^{-1}}=\displaystyle4\pi\int_0^1\frac{1}{r^{1/2}}dr=8\pi,\;\;\;
\|f\|_{\dot H^{1/2}}^2=\displaystyle4\pi\int_0^1\frac{1}{r}dr=\infty$$
and
$$\|g\|_{\bf {\mathcal X}^{-1}}=\displaystyle4\pi\int_1^\infty\frac{1}{r^{3/4}}dr=\infty,\;\;\;\|g\|_{\dot H^{1/2}}^2=\displaystyle4\pi\int_1^\infty\frac{1}{r^{3/2}}dr=8\pi.$$
\end{rem}
\section{Proof of Theorem \ref{jamel1}} This proof is inspired from the work of Gallagher-Iftimie-Planchon in \cite{GIP}.\\
Let $\varepsilon>0$, a sufficient condition on $\varepsilon$ is as follows
\begin{equation}\label{eqv1}\varepsilon\leq \frac{\nu}{2}.\end{equation}
For $k\in\mathbb N$, put
$${\mathcal A}_k=\{\xi\in\mathbb R^3;\;\;\;|\xi|\leq k\;\;\;{\rm and}\;\;\;|\widehat{u^0}(\xi)|\leq k\}$$
Clearly $\mathcal F^{-1}({\bf 1}_{{\mathcal A}_k}\widehat{u^0})$ converges to $u^0$ in ${\mathcal X}^{-1}(\mathbb R^3)$. Then, there is $k\in\mathbb N$ such that
$$\|u^0-\mathcal F^{-1}({\bf 1}_{{\mathcal A}_k}\widehat{u^0})\|_{{\mathcal X}^{-1}}<\varepsilon/2.$$
Put $v_k^0$ and $w_k^0$ as follows
$$\begin{array}{ccc}
  v_k^0&=& \mathcal F^{-1}({\bf 1}_{{\mathcal A}_k}\widehat{u^0})\\
 w_k^0&=& u^0- v_k^0.
  \end{array}
$$
Then $\| w_k^0\|_{{\mathcal X}^{-1}}<\varepsilon/2$ and $v_k^0\in{\mathcal X}^{-1}(\mathbb R^3)\cap L^2(\mathbb R^3) $. Now, consider the following system
$$ \left\{\begin{array}{l}
  \displaystyle\partial_tw-  \nu\Delta w+w.\nabla w=-\nabla p_{1,k},\quad
\mbox{in}\;\; \mathbb{R}^+\times\mathbb R^3\\
\mbox{div}\;   w= 0 \quad
\mbox{in}\; \; \mathbb{R}^+\times\mathbb R^3  \\
w _{\mid t=0}=w_k^0 \quad\mbox{in}\;\;\mathbb R^3.
\end{array}\right.\leqno(NS_k)$$
As $\| w_k^0\|_{{\mathcal X}^{-1}}<\varepsilon/2<\nu$ and by using Theorem \ref{zhen} and inequality (\ref{jam1}), we get a unique global solution $w_k$ of $(NS_k)$ such that $w_k\in{\mathcal C}(\mathbb R^+,{\bf {\mathcal X}^{-1}}\mathbb R^3)\cap L^1(\mathbb R^+,{\bf {\mathcal X}^{1}}\mathbb R^3).$ Moreover,
\begin{equation}\label{smid1}\|w_k(t)\|_{\bf {\mathcal X}^{-1}}+\frac{\nu}{2}\int_0^t\|w_k(z)\|_{\bf {\mathcal X}^1}dz\leq\|w^0_k\|_{\bf {\mathcal X}^{-1}},\,\forall t\geq 0.\end{equation}
Put $v_k=u-w_k$, clearly $v_k\in{\mathcal C}(\mathbb R^+,{\bf {\mathcal X}^{-1}}(\mathbb R^3))$ and satisfies
$$ \left\{\begin{array}{l}
  \displaystyle\partial_t
v_k-\nu\Delta v_k+ v_k.\nabla v_k+ w_k.\nabla v_k+ v_k.\nabla w_k=-\nabla p_{2,k},\quad
\mbox{in}\;\; \mathbb{R}^+\times\mathbb R^3\\
\mbox{div}\,v_k = 0 \quad
\mbox{in}\; \; \mathbb{R}^+\times\mathbb R^3  \\
{v_k}_{\mid t=0}=v^0_k \quad\mbox{in}\;\;\mathbb R^3.
\end{array}\right.$$
Taking the inner product in $L^2(\mathbb R^3)$ with $v_k$, we get
$$\frac{1}{2}\frac{d}{dt}\|v_k\|_{L^2}^2+\nu\|\nabla v_k\|_{L^2}^2\leq |\langle w_k.\nabla v_k/v_k\rangle_{L^2}|.$$
To estimate the RHT,
$$\begin{array}{ccc}
|\langle w_k.\nabla v_k/v_k\rangle_{L^2}|&=&|\langle{\rm div}\,( w_k\otimes v_k)/v_k\rangle_{L^2}|\\
&=&|\langle w_k\otimes v_k/\nabla v_k\rangle_{L^2}|\\
&\leq&\| w_k\otimes v_k\|_{L^2}\|\nabla v_k\|_{L^2}|\\
&\leq&\| \mathcal F(w_k\otimes v_k)\|_{L^2}\|\nabla v_k\|_{L^2}\\
&\leq&\| \widehat{w_k}*\widehat{v_k}\|_{L^2}\|\nabla v_k\|_{L^2}.
  \end{array}
$$
Young inequality and Lemma \ref{scl1} give
$$\begin{array}{ccc}
|\langle w_k.\nabla v_k/v_k\rangle_{L^2}|&\leq&\| \widehat{w_k}\|_{L^1}\|\widehat{v_k}\|_{L^2}\|\nabla v_k\|_{L^2}\;\;\;\;\;\;\;\;\\
 &\leq&\| w_k\|_{\bf {\mathcal X}^{-1}}^{1/2}\| w_k\|_{\bf {\mathcal X}^{1}}^{1/2}\|v_k\|_{L^2}\|\nabla v_k\|_{L^2}.
  \end{array}
$$
Using inequality $ab\leq \frac{a^2}{2}+\frac{b^2}{2}$, we get
$$|\langle w_k.\nabla v_k/v_k\rangle_{L^2}|\leq \frac{1}{2\nu}\| w_k\|_{\bf {\mathcal X}^{-1}}\| w_k\|_{\bf {\mathcal X}^{1}}\|v_k\|_{L^2}^2+\frac{\nu}{2}\|\nabla v_k\|_{L^2}$$
and
$$\frac{d}{dt}\|v_k\|_{L^2}^2+\nu\|\nabla v_k\|_{L^2}^2\leq\frac{1}{\nu}\| w_k\|_{\bf {\mathcal X}^{-1}}\| w_k\|_{\bf {\mathcal X}^{1}}\|v_k\|_{L^2}^2.$$
Gronwall Lemma yields
$$\|v_k\|_{L^2}^2+\nu\int_0^t\|\nabla v_k\|_{L^2}^2\leq\|v_k^0\|_{L^2}^2e^{\frac{1}{\nu}\int_0^t\| w_k\|_{\bf {\mathcal X}^{-1}}\| w_k\|_{\bf {\mathcal X}^{1}}}.$$
Using inequality (\ref{smid1}) , we get
$$\|v_k\|_{L^2}^2+\nu\int_0^t\|\nabla v_k\|_{L^2}^2\leq\|v_k^0\|_{L^2}^2e^{\frac{2}{\nu^2}\| w_k^0\|_{\bf {\mathcal X}^{-1}}^2}.$$
Combining the above inequality and Lemma \ref{scl2}, we can deduce that $v_k\in L^4(\mathbb R^+,{\bf{\mathcal X}^{-1}}(\mathbb R^3))$, and
$$\int_0^\infty\| v_k\|_{\bf {\mathcal X}^{-1}}^4\leq \int_0^\infty \|v_k\|_{L^2}^2\|\nabla v_k\|_{L^2}^2\leq \nu^{-1}\|v_k^0\|_{L^2}^4e^{\frac{4}{\nu^2}\| w_k^0\|_{\bf {\mathcal X}^{-1}}^2}.$$
By continuity of $v_k$ in ${\bf{\mathcal X}^{-1}}(\mathbb R^3)$, there is a time $t_0$ such that $\| v_k(t_0)\|_{\bf {\mathcal X}^{-1}}<\varepsilon/2$. Using equation (\ref{smid1}), we get
$$\| u(t_0)\|_{\bf {\mathcal X}^{-1}}\leq\| v_k(t_0)\|_{\bf {\mathcal X}^{-1}}+\| w_k(t_0)\|_{\bf {\mathcal X}^{-1}}<\frac{\varepsilon}{2}+\frac{\varepsilon}{2}=\varepsilon.$$
Now, consider the Navier-Stokes  system  starting at $t=t_0$ and using the global existence for the small initial data, we get
$$\|u(t)\|_{\bf {\mathcal X}^{-1}}+(\nu-\varepsilon)\int_{t_0}^t\|\Delta
u(\tau)\|_{\bf {\mathcal X}^{-1}}d\tau\leq  \varepsilon,\;\forall t\geq t_0.$$
Then, the desired result is proved.
\section{Stability of global solutions}
In this section we prove Theorem \ref{jamel2}.  This proof is done in two steps.\\
 {\bf Step 1:} Beginning by proving the following property: If $u$ is a maximal solution of $(NS)$ system with $u^0\in{\bf {\mathcal X}^{-1}}(\mathbb R^3)$ and $T^*$ is the maximal time of existence. We know that $u\in{\mathcal C}([0,T^*);{\bf {\mathcal X}^{-1}}(\mathbb R^3))\cap L^1_{loc}([0,T^*),{\bf {\mathcal X}^1}(\mathbb R^3))$. We have, if $T^*<\infty$ then
$$\int_0^{T^*}\|u(t)\|_{{\bf {\mathcal X}^1}}dt=\infty.$$
Indeed: Suppose that $\int_0^{T^*}\|u(t)\|_{{\bf {\mathcal X}^1}}dt<\infty$. Let a time $T\in(0,T^*)$ such that $\int_T^{T^*}\|u(t)\|_{{\bf {\mathcal X}^1}}dt<1/2$. Lemma \ref{scl1} gives, for all $t\in[T,T^*)$,
$$\begin{array}{ccc}
\displaystyle\|u(t)\|_{{\bf {\mathcal X}^{-1}}}+\int_T^t\|u(s)\|_{{\bf {\mathcal X}^1}}ds
&\leq&\displaystyle\|u(T)\|_{{\bf {\mathcal X}^{-1}}}+\int_T^t\|u(s)\|_{{\bf {\mathcal X}^{-1}}}\|u(s)\|_{{\bf {\mathcal X}^1}}ds\;\;\;\;\;\;\;\;\;\;\;\\
&\leq&\displaystyle\|u(T)\|_{{\bf {\mathcal X}^{-1}}}+\sup_{z\in[T,t]}\|u(z)\|_{{\bf {\mathcal X}^{-1}}}\int_T^t\|u(s)\|_{{\bf {\mathcal X}^1}}ds\\
&\leq&\displaystyle\|u(T)\|_{{\bf {\mathcal X}^{-1}}}+\frac{1}{2}\sup_{z\in[T,t]}\|u(z)\|_{{\bf {\mathcal X}^{-1}}}.\;\;\;\;\;\;\;\;\;\;\;\;\;\;\;\;\;\;\;\;\\
  \end{array}
$$
We can deduce $$\|u(s)\|_{{\bf {\mathcal X}^{-1}}}\leq 2\|u(T)\|_{{\bf {\mathcal X}^{-1}}},\;\;\forall s\in[T,T^*).$$
Let $M=\sup_{z\in[T,T^*)}\|u(z)\|_{{\bf {\mathcal X}^{-1}}}<\infty.$ We have
$$u(t')-u(t)=\nu\int_t^{t'}\Delta u-\int_t^{t'}{\mathbb P}{\rm div}\,(u\otimes u)$$
Using Lemma \ref{scl1}, we get
$$\begin{array}{ccc}
\displaystyle \|u(t')-u(t)\|_{{\bf {\mathcal X}^{-1}}}&\leq&\displaystyle\nu\int_t^{t'}\|u(s)\|_{{\bf {\mathcal X}^1}}ds+\int_t^{t'}\|u(s)\|_{{\bf {\mathcal X}^{-1}}}\|u(s)\|_{{\bf {\mathcal X}^1}}ds\\
&\leq&\displaystyle(\nu+M)\int_t^{t'}\|u(s)\|_{{\bf {\mathcal X}^1}}ds\;\;\;\;\;\;\;\;\;\;\;\;\;\;\;\;\;\;\;\;\;\;\;\;\;\;\;\;\;\;\;\;\;\;\;\;\;\\
\end{array}
$$
where the RHT goes to zero as $t$ and $t'$ tends to $T^*$. Then $u(t)$ is a Cauchy type at $T^*$. As ${\bf {\mathcal X}^{-1}}(\mathbb R^3)$ is Banach space, then there is an element $u^*$ in ${\bf {\mathcal X}^{-1}}(\mathbb R^3)$ such that $u(t)\rightarrow u^*$ in ${\bf {\mathcal X}^{-1}}(\mathbb R^3)$ if $t$ goes to $T^*$. Now, consider the Navier-Stokes system starting by $u^*$, using Theorem \ref{ZYth}, we get a unique solution which is extends $u$ beyond to $T^*$ which is absurd.\\\\
{\bf Step 2:} Let $v\in {\mathcal C}([0,T^*), {\bf {\mathcal X}^{-1}}(\mathbb R^3))$ be the maximal solution of $(NS)$ corresponding to the initial condition $v^0$. We want to prove $T^*=\infty$. Beginning by using Theorem \ref{ZYth}, we get  $v\in L^1_{loc}([0,T^*),{\bf {\mathcal X}^{1}}(\mathbb R^3))$.\\
 Put $w=v-u$ and $w^0=v^0-u(0)$. We have
 $$\partial_t w-\nu \Delta w+w.\nabla w+u.\nabla w+w.\nabla u=-\nabla P$$
 or
  $$\partial_t w-\nu \Delta w+{\rm div}\,(w\otimes w)+{\rm div}\,(u\otimes w)+{\rm div}\,(w\otimes u)=-\nabla P.$$
  Then, for $t\in[0,T^*)$
  $$\|w(t)\|_{{\bf {\mathcal X}^{-1}}}+\nu\int_0^t\|w(t)\|_{{\bf {\mathcal X}^{1}}}\leq \|w^0\|_{{\bf {\mathcal X}^{-1}}}+(I)+(II)$$
  where
  $$\begin{array}{ccc}
         (I)&=&\displaystyle\int_0^t\|{\rm div}\,(w\otimes w)\|_{{\bf {\mathcal X}^{-1}}} \\
         (II)&=&\displaystyle\int_0^t\|{\rm div}\,(u\otimes w)\|_{{\bf {\mathcal X}^{-1}}}+\|{\rm div}\,(w\otimes u)\|_{{\bf {\mathcal X}^{-1}}}.
       \end{array}
  $$
Lemma \ref{scl1} gives
 $$\begin{array}{ccc}
   (I)&\leq &\displaystyle\int_0^t\|w\otimes w\|_{{{\bf {\mathcal X}^0}}} \\
   &\leq &\displaystyle\int_0^t\|w\|_{{\bf {\mathcal X}^{-1}}}\|w\|_{{\bf {\mathcal X}^{1}}},
   \end{array}
 $$
 and
$$\begin{array}{ccc}
   (II)&\leq &\displaystyle\int_0^t\|u\otimes w\|_{{\bf {\mathcal X}^0}}+\|w\otimes u\|_{{\bf {\mathcal X}^0}} \\
   &\leq &2\displaystyle\int_0^t\|u\|_{{\bf {\mathcal X}^0}}\|w\|_{{\bf {\mathcal X}^{-1}}}^{1/2}\|w\|_{{\bf {\mathcal X}^{1}}}^{1/2}\\
    &\leq &\displaystyle\frac{4}{\nu}\int_0^t\|u\|_{{\bf {\mathcal X}^0}}^2\|w\|_{{\bf {\mathcal X}^{-1}}}+\frac{\nu}{4}\int_0^t\|w\|_{{\bf {\mathcal X}^{1}}}.
   \end{array}
 $$
 Then
$$\|w(t)\|_{{\bf {\mathcal X}^{-1}}}+\frac{3\nu}{4}\int_0^t\|w(t)\|_{{\bf {\mathcal X}^{1}}}\leq \|w^0\|_{{\bf {\mathcal X}^{-1}}}+\int_0^t\|w\|_{{\bf {\mathcal X}^{-1}}}\|w\|_{{\bf {\mathcal X}^{1}}}+\frac{2}{\nu}\int_0^t\|u\|_{{\bf {\mathcal X}^0}}^2\|w\|_{{\bf {\mathcal X}^{-1}}}.$$
Put $$T=\sup\{t\in[0,T^*),\;\sup_{z\in[0,t]}\|w(z)\|_{{\bf {\mathcal X}^{-1}}}<\frac{\nu}{4}\}.$$
For $t\in[0,T)$, we have
$$\|w(t)\|_{{\bf {\mathcal X}^{-1}}}+\frac{\nu}{2}\int_0^t\|w(t)\|_{{\bf {\mathcal X}^{1}}}\leq \|w^0\|_{{\bf {\mathcal X}^{-1}}}+\frac{2}{\nu}\int_0^t\|u\|_{{\bf {\mathcal X}^0}}^2\|w\|_{{\bf {\mathcal X}^{-1}}}.$$
Gronwall Lemma yields
$$\|w(t)\|_{{\bf {\mathcal X}^{-1}}}+\frac{\nu}{2}\int_0^t\|w(t)\|_{{\bf {\mathcal X}^{1}}}\leq \|w^0\|_{{\bf {\mathcal X}^{-1}}}e^{\frac{2}{\nu}\int_0^t\|\widehat{u}\|_{L^1}^2}\leq \|w^0\|_{{\bf {\mathcal X}^{-1}}}e^{\frac{2}{\nu}\int_0^\infty\|\widehat{u}\|_{L^1}^2}<\frac{\nu}{8}.$$
Then $T=T^*$ and $\int_0^{T^*}\|w(t)\|_{{\bf {\mathcal X}^{1}}}<\infty$, therefore $T^*=\infty$ and the proof is finished.

\end{document}